\documentclass[12pt]{amsart}%\documentclass[12pt,a4paper]{amsart}
\usepackage{amsfonts}

\usepackage[letterpaper]{geometry}
\usepackage{xcolor}
\usepackage{hyperref}

\usepackage{amssymb}
\def\R{\mathbb{R}}

\def\N{\mathbb{N}}

\newtheorem{thm}{\bf Theorem}[section]

\newtheorem{prop}[thm]{\bf Proposition}

\theoremstyle{definition}
\newtheorem*{defi}{\bf Definition}
\newtheorem*{rem}{\bf Remark}

\begin{document}

\title[Consecutive sums and Reconstructing Problems]{Consecutive sums and Reconstructing Problems}

 \author[Norbert Hegyv\'ari]{Norbert Hegyv\'ari}
 \address{Norbert Hegyv\'{a}ri, ELTE TTK,
E\"otv\"os University, Institute of Mathematics, H-1117
P\'{a}zm\'{a}ny st. 1/c, Budapest, Hungary and associated member of Alfr\'ed R\'enyi Institute of Mathematics, Hungarian Academy of Science, H-1364 Budapest, P.O.Box 127.}
 \email{hegyvari@renyi.hu}

\author[P\'eter P\'al Pach]{P\'eter P\'al Pach}
\address{P\'eter P\'al Pach, HUN-REN Alfr\'ed R\'enyi Institute of Mathematics, Re\'altanoda utca 13--15., H-1053 Budapest,  Hungary; \newline \hspace*{4mm}
MTA--HUN-REN RI Lend\"ulet ``Momentum'' Arithmetic Combinatorics Research Group, Re\'altanoda utca 13--15., H-1053 Budapest,  Hungary; \newline \hspace*{4mm}
Department of Computer Science and Information Theory, Budapest University of Technology and Economics, M\H{u}egyetem rkp. 3., H-1111 Budapest, Hungary.}
 \email{pachpp@renyi.hu}

\author[Quanyu Tang]{Quanyu Tang}
\address{Quanyu Tang, School of Mathematical Sciences, University of Science and Technology of China, Hefei 230026, P. R. China}
 \email{tangquanyu827@gmail.com}

\begin{abstract}
    We study several problems concerning consecutive sums of integer sequences and their relationship with Sidon sets. We also provide new lower bounds for an extremal function related to sequences with distinct consecutive sums, and consider reconstruction problems from multisets of consecutive sums.
\end{abstract}

\maketitle

\section{Introduction}

A classical problem, posed by Erdős and Harzheim, asks whether one can find a sequence of integers bounded by $n$ such that all of its consecutive sums are pairwise different. The first named author addressed this problem by revealing a connection to Sidon sets: the condition that all consecutive sums are distinct implies that the sequence of their prefix sums forms a Sidon sequence.

These results, together with subsequent developments in inverse additive combinatorics (cf. for instance Beker \cite{B}, Tao \cite{T}), demonstrate how structural properties of a sequence can be deduced from information about its consecutive sums.

In this paper, we investigate several related questions: ranging from connections between convexity and Sidon properties, extremal problems and reconstruction algorithms (see related results [K], [N]).

The paper is organized as follows. In Section~2, we investigate convex sequences and their relationship with Sidon sets. In Section~3, we establish lower bounds for the extremal function $f(a,K)$ related to consecutive sums. Finally, in Section~4, we study reconstruction from consecutive sums in one dimension and from block-sums in two dimensions.

\section{Convex sequences vs. Sidon sets}

Throughout the paper, we use the convention \(\mathbb N=\{1,2,3,\ldots\}\). We use the standard definition of a Sidon set: a set
\(A\subset\mathbb Z\) is called a Sidon set if
\[
x_1+x_2=x_3+x_4,
\qquad x_1,x_2,x_3,x_4\in A,
\]
implies \(\{x_1,x_2\}=\{x_3,x_4\}\) as unordered multisets. Repeated summands are allowed.

A classical example of a convex sequence is the sequence of squares, which is not a Sidon sequence.
Of course, there are also sequences that are both convex and Sidon.  The simplest example is $T=\{2^n\}_{n\geq 0}$.
 A denser sequence is due to Ruzsa.
He proved that there exists an irrational number $c$ (where $0 < c < 1$) such that the sequence generated by the following formula forms an infinite Sidon sequence:
$$r_n = n^5 + \lfloor c \cdot n^4 \rfloor.$$
In particular, there exists a fixed $N_0$ such that $\{r_n\}_{n\geq N_0}$ is a convex sequence. More generally, he proved the following theorem.

\begin{thm}[Ruzsa \cite{R}]

Let $\alpha, \beta$ be real numbers satisfying $ \alpha > 4, \ \alpha > \beta > 3 + 1/(\alpha - 1).$

Write $f(n) = \lfloor n^\alpha + \xi n^\beta\rfloor.$ For almost every $\xi \in [0, 1]$ there exists some $n_0$ (depending on $\alpha, \beta, \xi$) such that the set $A = \{f(n) : n > n_0\}$
is a Sidon set.
\end{thm}

The question arises: What is the degree of a generalized polynomial that is a Sidon sequence?

\begin{defi}
    Generalized polynomials are functions (or sequences) constructed from conventional polynomials, with the allowance of an additional operation: the application of the floor function ($\lfloor \dots \rfloor$) or the fractional part function ($\{\dots\}$).
\end{defi}
For example $\lfloor \sqrt{2} x \rfloor, \  \sqrt{3} x \lfloor x^2 + \lfloor \sqrt{5} x \rfloor \rfloor$ are generalized polynomials. In the rest we consider generalized polynomials $P: \N \to \R$.

We now study a class of quadratic generalized polynomials for which
strict convexity can be guaranteed, and ask whether they can form
infinite Sidon sequences.

\begin{defi}
Let $\mathcal P$ be the class of integer-valued sequences of the form
\[
P(n)=c_1n^2+\sum_{i=1}^{k}c_{2,i}
\left\lfloor\alpha_i n^2+\nu_i n+\mu_i\right\rfloor+L(n),
\qquad n\in\mathbb N,
\]
where $k\ge0$ is fixed, all coefficients are real, and $L(n)$ is a
generalized polynomial of degree at most one satisfying
\[
D_L:=\sup_{n\ge1}\bigl|L(n+2)-2L(n+1)+L(n)\bigr|<\infty.
\]
We additionally require
\[
c_1+\sum_{i=1}^{k}c_{2,i}\alpha_i
>\sum_{i=1}^{k}|c_{2,i}|+\frac{D_L}{2}.
\]
\end{defi}

\begin{thm}
Every sequence $P\in\mathcal P$ is strictly convex.
Moreover, the set \(\{P(n):n\geq1\}\) is not a Sidon set.
\end{thm}
\begin{proof}
Write
\[
c=c_1+\sum_{i=1}^{k}c_{2,i}\alpha_i,
\qquad q_i(n)=\alpha_i n^2+\nu_i n+\mu_i,
\]
and let $\Delta s(n)=s(n+1)-s(n)$.
For $\theta_i(n)=\{q_i(n)\}\in[0,1)$, we have
$|\Delta^2\theta_i(n)|<2$ and \(\Delta^2\lfloor q_i(n)\rfloor=2\alpha_i-\Delta^2\theta_i(n)\). Consequently,
\[
\Delta^2P(n)=2c-\sum_{i=1}^{k}c_{2,i}
\Delta^2\theta_i(n)+\Delta^2L(n).
\]
Set
\[
\delta=2c-2\sum_{i=1}^{k}|c_{2,i}|-D_L>0,
\qquad
B=2c+2\sum_{i=1}^{k}|c_{2,i}|+D_L.
\]
Then for every $n\ge1$,
\[
0<\delta\le\Delta^2P(n)\le B<\infty.
\]
Thus $P$ is strictly convex. In particular, its first differences
$\Delta P(n)$ are strictly increasing and tend to infinity, so
there exists $n_0$ such that $\Delta P(n)>0$ for all $n\ge n_0$.

For $r\ge1$, let
\[
b_r=P(n_0+r)-P(n_0+r-1).
\]
The $b_r$ are positive integers, and
\[
0<b_{r+1}-b_r=\Delta^2P(n_0+r-1)\le B.
\]
Choose a positive integer $K\ge B$.
If the values of $P$ formed an infinite Sidon set, then all the
consecutive sums of the sequence $(b_r)_{r\ge1}$ would be distinct:
indeed, these sums are the positive differences
\[
\sum_{r=u}^{v}b_r=P(n_0+v)-P(n_0+u-1),
\qquad 1\le u\le v,
\]
and a Sidon set has distinct positive differences.
Therefore, by Theorem~\ref{thm:f-upper}, we have
$b_m\le f(b_1,K)$ for every $m\ge1$.
But $b_m\ge b_1+m-1$ because the $b_r$ are strictly increasing
integers. This contradiction proves that $P$ cannot form an
infinite Sidon set.
\end{proof}

\begin{rem}
    One might ask why we set a condition on the coefficients in the definition of $\mathcal P$. To see why some
restriction is needed, observe that for every $k\in\mathbb N$
and every sufficiently large integer $n$ (for instance, $n>8k$),
there exists $\alpha>0$ such that
\[
\lfloor\alpha n^2\rfloor=k,\qquad
\lfloor\alpha(n+1)^2\rfloor=k+1,\qquad
\lfloor\alpha(n+2)^2\rfloor=k+1.
\]

Let $\alpha = \frac{k+1}{(n+1)^2} + \frac{1}{(n+2)^4}$.

We claim
$$k \le \frac{k+1}{(n+1)^2}n^2 + \frac{1}{(n+2)^4} n^2 < k+1.$$
For the lower bound we need $n^2 \ge 2kn + k$; the upper bound by a short calculation holds.
Hence $\left\lfloor \alpha n^2\right\rfloor=k$. For $\left\lfloor \alpha (n+1)^2\right\rfloor=k+1$ we need
$$k+1 \le (k+1) + \frac{1}{(n+2)^4}(n+1)^2 < k+2.$$
Both the lower and upper bounds are trivial.
Finally, for $k+1 \le \alpha(n+2)^2 < k+2$ we need
$$k+1 \le \frac{k+1}{(n+1)^2}(n+2)^2 +  \frac{1}{(n+2)^4}(n+2)^2 < k+2,$$
the lower bound holds. For the upper bound we need
$\frac{1}{(n+2)^4}(n+2)^2 < (k+2) - \frac{k+1}{(n+1)^2}(n+2)^2$, from this
$$
\frac{1}{(n+2)^2}< \frac{n^2-2kn-3k-2}{(n+1)^2}
$$
which holds when $n^2-2kn-3k>3$.
From the previous requirement namely that $n^2 \ge 2kn + k$, if $n>8k$, then these inequalities hold.
\end{rem}

\section{Lower bound for $f(a,K)$}

In \cite{H} the first named author answered the following question of Erd\H os:

Question: Let $a_1 < a_2 < \dots < a_n$ be a strictly increasing sequence of positive integers with initial term $a_1 = a$. Suppose that the difference between consecutive terms is bounded by a positive integer $K$, that is, $a_{i+1} - a_i \le K$ for all $1 \le i \le n-1$. If all consecutive sums $\sum_{i=u}^{v} a_i$ (where $1 \le u \le v \le n$) are distinct, what is the maximum possible value of the last element $a_n$?

The first named author~\cite[Theorem~3]{H} proved the following upper bound.
\begin{thm}[Hegyv\'ari \cite{H}]\label{thm:f-upper}
Let \(a,K\geq 1\) be integers. Let $f(a,K)$ be the maximum integer $a_n$ for which there exists a sequence $a = a_1 < a_2 < \dots < a_n$ such that $a_{i+1} - a_i \le K$ ($1 \le i \le n-1$) and all consecutive sums are different. Then$$f(a,K) < \left(a + \frac{K}{2}\right) e^{K+1} + K e^{2K+2}.$$
\end{thm}
Note that a more careful calculus shows the bound
$f(a,K) < 2e^{K+0.5}(a+K/2) + \frac{K}{2}e^{2K+1}$.

Here we provide lower bounds for $f(a,K)$.

\begin{thm}
Let $a, K \ge 1$ be integers with $a \ge K$. Then
\[ f(a,K) \ge a + K \left( \left\lfloor 2\sqrt{\frac{a}{K}-1} \right\rfloor - 1 \right) \]
In particular, as $a/K \rightarrow \infty$,
\[ f(a,K) \ge a + (2-o(1))\sqrt{aK}. \]
\end{thm}

\begin{proof}
Set
\[ k := \left\lfloor 2\sqrt{\frac{a}{K}-1} \right\rfloor. \]
If $k=0$, then the stated lower bound is trivial, so we may assume $k \ge 1$.
Consider the arithmetic progression
\[ a_i := a + (i-1)K \quad (1 \le i \le k). \]
Then $a_1 = a$, $a_{i+1} - a_i = K$ $(1 \le i \le k-1)$, and $a_k = a + (k-1)K$.
Therefore it is enough to prove that all consecutive sums of $a_1, \dots, a_k$ are distinct.

Write $a_i = K(t+i)$, where $t := \frac{a}{K}-1$. By the choice of $k$, we have $t \ge \frac{k^2}{4}$.

Suppose, for contradiction, that two distinct consecutive sums are equal. Then there exist integers $u, v, r, s$ with
\[ 0 \le u < u+r \le k, \quad 0 \le v < v+s \le k, \quad (u,r) \ne (v,s), \]
such that
\[ \sum_{i=u+1}^{u+r} a_i = \sum_{j=v+1}^{v+s} a_j. \]
We first note that $r \ne s$. Indeed, for any fixed length $\ell$, the sums $a_{m+1} + \dots + a_{m+\ell}$ $(0 \le m \le k-\ell)$ form a strictly increasing sequence in $m$.%, because

Thus, after interchanging the two blocks if necessary, we may assume that $r < s$. Dividing the equality of sums by $K$, we obtain
\[ \sum_{i=u+1}^{u+r} (t+i) = \sum_{j=v+1}^{v+s} (t+j), \]
so
\[ \sum_{i=u+1}^{u+r} i - \sum_{j=v+1}^{v+s} j = (s-r)t. \]
Since $s-r \ge 1$, the right-hand side is at least $t$.

On the other hand, among all choices of a block of length $r$ and a block of length $s$, the quantity
\[ \sum_{i=u+1}^{u+r} i - \sum_{j=v+1}^{v+s} j \]
is maximal when the $r$-block is as far to the right as possible and the $s$-block is as far to the left as possible. Therefore
\[ \sum_{i=u+1}^{u+r} i - \sum_{j=v+1}^{v+s} j \le \sum_{i=k-r+1}^{k} i - \sum_{j=1}^{s} j. \]
Hence
\[ \sum_{i=u+1}^{u+r} i - \sum_{j=v+1}^{v+s} j \le \frac{r(2k-r+1)-s(s+1)}{2}. \]
Since $s \ge r+1$, it follows that
\[ \frac{r(2k-r+1)-s(s+1)}{2} \le \frac{r(2k-r+1)-(r+1)(r+2)}{2} = r(k-r-1)-1. \]
Now $r(k-r-1)  \le \frac{(k-1)^2}{4}$, and consequently
\[ \sum_{i=u+1}^{u+r} i - \sum_{j=v+1}^{v+s} j < \frac{k^2}{4} \le t. \]
This contradicts
\[ \sum_{i=u+1}^{u+r} i - \sum_{j=v+1}^{v+s} j = (s-r)t \ge t. \]
Therefore, all consecutive sums of $a_1, \dots, a_k$ are distinct. Since this sequence starts at $a$, has all gaps equal to $K$, and ends at $a_k$, the definition of $f(a,K)$ gives
\[ f(a,K) \ge a_k = a + (k-1)K = a + K\left( \left\lfloor 2\sqrt{\frac{a}{K}-1} \right\rfloor - 1 \right). \]
Finally, if $a/K \rightarrow \infty$, then
\[ \left\lfloor 2\sqrt{\frac{a}{K}-1} \right\rfloor - 1 = 2\sqrt{\frac{a}{K}} - o\left(\sqrt{\frac{a}{K}}\right), \]
and hence $f(a,K) \ge a + (2-o(1))\sqrt{aK}$. This completes the proof.
\end{proof}

We give another lower bound which is better when $a\ll K$.

\begin{thm}
   If $K$ is sufficiently large, then $f(a,K)\geq a+\frac14K^2$.
\end{thm}

\begin{proof}
   Assume first that there is a prime $p\in(K/2,K]$ such that $p\nmid a$. Let $a_i=a+(i-1)p$ for $1\leq i\leq n:=p+1$. Each $i$-term consecutive sum has residue $ia+\binom{i}{2}p$ modulo $ip$. Clearly, all the consecutive $i$-term sums are distinct. If an $i$-term sum is equal to a $j$-term sum (where $i<j$), then $p\mid ia-ja=(i-j)a$. Since $p\nmid a$, this is not possible. (Indeed, $i=1,j=p+1$ is not possible, since $a_1+\dots+a_{p+1}$ is larger than all 1-term sums.) Therefore, in this case $f(a,K)\geq a+p^2> a+\frac{1}{4}K^2$.

   Now, assume that every prime $p\in(K/2,K]$ divides $a$. Let $a_1=a$ and $a_i=a+(i-1)p-1$ for $1<i\leq n:=p+1$. Each $i$-term consecutive sum has residue $i(a-1)+\binom{i}{2}p$ or $i(a-1)+\binom{i}{2}p+1$  modulo $ip$. Let $i<j$ and assume that an $i$-term consecutive sum is equal to a $j$-term consecutive sum.    It is not possible that the $i$-term sum contains $a_1$, since this would be the smallest $i$-term sum, smaller than all the $j$-term sums.

   If none of the two sums contain $a_1$, then $p\mid (i-j)(a-1)$, which is not possible, since $p\mid a$.

   The only remaining case is when the $j$-term sum contains $a_1$ (and the $i$-term sum does not contain $a_1$). This leads to the congruence $i(a-1)\equiv j(a-1)+1\pmod{p}$, implying $i\equiv j-1 \pmod{p}$, thus $i=j-1$.

Let the $i$-term consecutive sum be $a_{s+1}+a_{s+2}+\dots+a_{s+i}$.
Then from the equation $a_1+\dots+a_{i+1}=a_{s+1}+a_{s+2}+\dots+a_{s+i}$ we get that $s=1+\frac{a}{pi}$. As $s+i\leq p+1$, we get that $a<ip^2<K^3$, however  the product of primes from $(K/2,K]$ divides the positive integer $a$ and for sufficiently large $K$ this product exceeds $K^3$. This contradiction implies that the consecutive sums are indeed distinct. Hence, in this case $f(a,K)\geq a+p^2-1>a+\frac{1}{4}K^2$ holds, as well.

\end{proof}

\section{Reconstruction in $\N$ and in $\N^2$}

The reconstruction of discrete structures from data is a fundamental problem in computer science, discrete mathematics, and inverse number theory. One classical example is discrete tomography, which focuses on reconstructing multidimensional grids, most commonly binary matrices. Similarly, in computational biology, the well-known DNA Reconstruction Problem asks for the reconstruction of a point set solely from the unordered multiset of its pairwise distances.

In this section, we study reconstruction problems related to consecutive sums in one and two dimensions.

\subsection{One-dimensional case}

Let $A=(1\leq a_1<a_2<\dots<a_n)$ be a strictly increasing sequence of integers. For
$1\le i<j\le n$, define the consecutive sum
\[
S_{i,j}:=\sum_{k=i}^{j}a_k.
\]
Let $C_1(A)$ denote the multiset of all consecutive sums, written in nondecreasing order as
\[
C_1(A)=\{x_1\le x_2\le\dots\le x_m\},
\qquad m=\binom{n}{2}.
\]
We consider only consecutive sums involving at least two terms.
\begin{prop}
    Let $n\geq 3$. Given $C_1(A)$, there is an algorithm with running time at most
    $O(n^2)$ that reconstructs the sequence $A$.
\end{prop}
\begin{proof}
Clearly, $x_m=S_{1,n}$ and $x_{m-1}=S_{2,n}$. Hence,
\[
a_1=S_{1,n}-S_{2,n}=x_m-x_{m-1}.
\]
Since the smallest element of $C_1(A)$ is $x_1=S_{1,2}=a_1+a_2$,
we obtain
\[
a_2=x_1-a_1.
\]

Now suppose that the elements $a_1<a_2<\dots<a_k$ have been
determined, and that we have calculated all consecutive sums
\[
\{S_{i,j}:1\leq i<j\leq k\}.
\]
Let $x_t$ be the smallest element of the multiset obtained from $C_1(A)$ by removing the previously determined consecutive sums, with their multiplicities. (Note that at each step, we remove the newly determined consecutive sums from the multiset, taking their multiplicities into account).  Since $a_1<\dots<a_n$, we have
$x_t=S_{k,k+1}=a_k+a_{k+1}.$

Consequently, $a_{k+1}=x_t-a_k.$
We can then calculate all the new consecutive sums
$S_{i,k+1}$ ($1\leq i\leq k$).

At the $(k+1)$-st step, we perform $k+1$ operations. Thus the total
number of operations is indeed $O(n^2)$.
\end{proof}

\subsection{Two-dimensional case}

In the one-dimensional case, we considered the sums over all intervals of consecutive indices of length at least two. We now introduce the corresponding two-dimensional notion.

\medskip

Let $X$ be an $N\times N$ matrix with positive real entries, $N\geq 3$. For every rectangular submatrix $X'$ of $X$ whose rows and columns are consecutive, let $S(X')$ denote the sum of its entries. We call these sums block-sums when the submatrix is not of size $1\times 1$.

In \cite{H}, the following translation problem was posed. Given a finite sequence $A$, determine the smallest integer $t$ such that all consecutive sums of $A+t$ are distinct. The existence of such an integer is immediate.

We consider the two-dimensional analogue. Let $M$ be an $N\times N$ matrix. Instead of adding an integer, we add an $N\times N$ translation matrix $T$.

\medskip

\textbf{The communication process.}
The specified reconstruction task involves the following communication process: given an $N\times N$ matrix $M$ with positive entries, the value of $S(M)$ is known. We send the translation matrix $T$,
and then receive the multiset of all block-sums of $M+T$. Our task is to reconstruct $M$ from this information.

\medskip

Let $\phi: \ [N]\times [N]\to \{1,2,\dots, N^2\}$ be any bijection (a fixed permutation of the entries of an $N\times N$ matrix) and let
$T_{ij}=K\cdot2^{\phi(i,j)}$, where $K>S(M)$.

\textbf{Recognizing the rectangles $T$.}

Let $T'$ be any rectangle in $T$. We have
$$
S(T')=K\sum_{(i,j)\in T'}2^{\phi(i,j)}.
$$
Since $K>S(M)$ and the binary representation of integers unique, from a sum $T'$ can be uniquely identified.

\smallskip

\textbf{Recognizing the rectangles of $M+T$}. Let $(M+T)'$ be any rectangle in $M+T$. Clearly $(M+T)'=M'+T'$ and $S((M+T)')=S(M')+S(T').$ Observe that
$$
S((M+T)')=S(M')+K\sum_{(i,j)\in T'}2^{\phi(i,j)},
$$
also, by $K>S(M)$ the contribution from $M$ is negligible, the sums determine which entries $(i,j)$ belong to  $T'$. Thus we can calculate each block-sum of $M$, as well (by subtrackting the corresponding block-sum of $T$ from the block-sum for $M+T$).

\medskip

First, $m_{11}$ is the difference of the two block-sums $m_{11}+m_{12}+m_{13}$ and $m_{12}+m_{13}$. Then, by using block-sums of the form $m_{ij}+m_{i,j+1}$ and $m_{ij}+m_{i+1,j}$ together with previously determined values we can calculate the entries of $M$ one by one.

The longest step is reconstructing the block-sums of $M+T$ which takes $O(N^4)$ steps.

Thus the following result is obtained.

\begin{thm}
    Let $N\geq 3$ and $M = (m_{i,j})$ be an $N \times N$ matrix such that $m_{i,j}>0$ for all $i, j$. Given the block-sums of the matrix $M+T$ for a suitably chosen and known translation matrix $T$, there is an algorithm that reconstructs $M$ within $O(N^4)$ steps.
\end{thm}

\section{Acknowledgements}
N.H. was supported by the National Research, Development and Innovation Office NKFIH Grant No K-146387. P.P.P. was supported by the National Research, Development and Innovation Office NKFIH (Grant Nr. K146387 and Excellence program, Grant Nr. 153829). Q.T. was supported by National Key Research and Development Program of China 2023YFA1010201.


\begin{thebibliography}{99}
\bibitem[B]{B} A. Beker, On a problem of Erdős and Graham about consecutive sums in strictly increasing sequences.  Bulletin of the London Mathematical Society 56(8) (2024) 2749--2759.

\bibitem[H]{H} N. Hegyv\'ari, On consecutive sums in sequences, Acta Math. Hungar. 48 (1986), 193--200.

\bibitem[K]{K} J. Konieczny, On consecutive sums in permutations, Journal of Combinatorics 12 (2021), no. 3, 413--477.

\bibitem[N]{N} M.B. Nathanson, Sidon sets with $\Delta-$separated sumsets in additive number theory, arXiv:2608.07416 [math.NT],
https://doi.org/10.48550/arXiv.2608.07416

\bibitem[R]{R} I. Z. Ruzsa, An almost polynomial Sidon sequence, Studia Scientiarum Mathematicarum Hungarica 38 (1-4) p. 367--376,
DOI:10.1556/SScMath.38.2001.1-4.27

\bibitem[T]{T} T. Tao,
 A set that represents all large integers multiple times by consecutive elements (2026)
 [Manuscript]. UCLA Department of Mathematics.
 Available at \url{https://terrytao.wordpress.com/wp-content/uploads/2026/02/erdos-358.pdf}
\end{thebibliography}
\end{document}